\documentclass[proceedings]{aofa}%

\usepackage[utf8]{inputenc}
\usepackage{subfigure}

\usepackage{amsmath,amssymb,amsthm,hyperref}
\usepackage{color}

\newtheorem{theorem}{Theorem}[section]

\theoremstyle{definition}

\author[M.~Kuba and A.~Panholzer]{Markus Kuba\addressmark{1} and Alois Panholzer\addressmark{2}\thanks{Supported by the Austrian Science Fund FWF, grant P25337-N23}}

\title[Analysis of growth models for series-parallel networks]{Combinatorial analysis of growth models for series-parallel networks}

\address{\addressmark{1}Institute of Applied Mathematics and Natural Sciences, University of Applied Sciences - Technikum Wien,
H\"och\-st\"adt\-platz 5, 1200 Wien, Austria, \email{kuba@technikum-wien.at}\\
\addressmark{2}Institut f{\"u}r Diskrete Mathematik und Geometrie, Technische Universit\"at Wien, Wiedner Hauptstr. 8-10/104, 1040 Wien, Austria, \email{Alois.Panholzer@tuwien.ac.at}}


\begin{document}

\maketitle

\begin{abstract}
We give combinatorial descriptions of two stochastic growth models for series-parallel networks introduced by Hosam Mahmoud by encoding the growth process via recursive tree structures. Using decompositions of the tree structures and applying analytic combinatorics methods allows a study of quantities in the corresponding series-parallel networks. For both models we obtain limiting distribution results for the degree of the poles and the length of a random source-to-sink path, and furthermore we get asymptotic results for the expected number of source-to-sink paths.
\end{abstract}

\section{Introduction}

Series-parallel networks are two-terminal graphs, i.e., they have two distinguished vertices called the source and the sink, that can be constructed recursively by applying two simple composition operations, namely the parallel composition (where the sources and the sinks of two series-parallel networks are merged) and the series composition (where the sink of one series-parallel network is merged with the source of another series-parallel network). Here we will always consider series-parallel networks as digraphs with edges oriented in direction from the north-pole, the source, towards the south-pole, the sink. Such graphs can be used to model the flow in a bipolar network, e.g., of current in an electric circuit or goods from the producer to a market. Furthermore series-parallel networks and series-parallel graphs (i.e., graphs which are series-parallel networks when some two of its vertices are regarded as source and sink; see, e.g., \cite{BraLeSpi1999} for exact definitions and alternative characterizations) are of interest in computational complexity theory, since some in general NP-complete graph problems are solvable in linear time on series-parallel graphs (e.g., finding a maximum independent set).

Recently there occurred several studies concerning the typical behaviour of structural quantities (as, e.g., node-degrees, see~\cite{DrmGimNoy2010}) in series-parallel graphs and networks under a uniform model of randomness, i.e., where all series-parallel graphs of a certain size (counted by the number of edges) are equally likely. In contrast to these uniform models, Mahmoud~\cite{Mahmoud2013,Mahmoud2014} introduced two interesting growth models for series-parallel networks, which are generated by starting with a single directed arc from the source to the sink and iteratively carrying out serial and parallel edge-duplications according to a stochastic growth rule; we call them uniform Bernoulli edge-duplication rule (``Bernoulli model'' for short) and uniform binary saturation edge-duplication rule (``binary model'' for short). A formal description of these models is given in Section~\ref{sec:Growth_models_recursive_trees}. Using the defining stochastic growth rules and a description via P\'{o}lya-Eggenberger urn models (see, e.g., \cite{Mah2009}), several quantities for series-parallel networks (as the number of nodes of small degree and the degree of the source for the Bernoulli model, and the length of a random source-to-sink path for the binary model) are treated in \cite{Mahmoud2013,Mahmoud2014}.

The aim of this work is to give an alternative description of these growth models for series-parallel networks by encoding the growth of them via recursive tree structures, to be precise, via edge-coloured recursive trees and so-called bucket-recursive trees (see \cite{KubPan2010} and references therein). The advantage of such a modelling is that these objects allow not only a stochastic description (the tree evolution process which reflects the growth rule of the series-parallel network), but also a combinatorial one (as certain increasingly labelled trees or bucket trees), which gives rise to a top-down decomposition of the structure. An important observation is that indeed various interesting quantities for series-parallel networks can be studied by considering certain parameters in the corresponding recursive tree model and making use of the combinatorial decomposition. We focus here on the quantities degree $D_{n}$ of the source and/or sink, length $L_{n}$ of a random source-to-sink path and the number $P_{n}$ of source-to-sink paths in a random series-parallel network of size $n$, but mention that also other quantities (as, e.g., the number of ancestors, node-degrees, or the number of paths through a random or the $j$-th edge) could be treated in a similar way. We obtain limiting distribution results for $D_{n}$ and $L_{n}$ (thus answering questions left open in \cite{Mahmoud2013,Mahmoud2014}), whereas for the r.v.\ $P_{n}$ (whose distributional treatment seems to be considerably more involved) we are able to give asymptotic results for the expectation. 

Mathematically, an analytic combinatorics treatment of the quantities of interest leads to studies of first and second order non-linear differential equations. In this context we want to mention that another model of series-parallel networks called increasing diamonds has been introduced recently in \cite{BodDieFonGenHWa2015+}. A treatment of quantities in such networks inherently also yields a study of second order non-linear differential equations; however, the definition as well as the structure of increasing diamonds is quite different from the models treated here as can be seen also by comparing the behaviour of typical graph parameters (e.g., the number of source-to-sink paths $P_{n}$ in increasing diamonds is trivially bounded by $n$, whereas in the models studied here the expected number of paths grows exponentially). We mention that the analysis of the structures considered here has further relations to other objects; e.g., it holds that the Mittag-Leffler limiting distributions occurring in Theorem~\ref{thm:Bernoulli_degree_limit} \& \ref{thm:Bernoulli_length_path} also appear in other combinatorial contexts as in certain triangular balanced urn models (see \cite{Janson2010}) or implicitly in the recent study of an extra clustering model for animal grouping \cite{DrmFucLee2015+} (after scaling, as continuous part of the characterization given in \cite[Theorem~2]{DrmFucLee2015+}, since it is possible to simplify some of the representations given there). Also the characterizations of the limiting distribution for $D_{n}$ and $L_{n}$ of binary series-parallel networks via the sequence of $r$-th integer moments satisfies a recurrence relation of ``convolution type'' similar to ones occurring in \cite{CheFerHwaMar2014}, for which asymptotic studies have been carried out.
Furthermore, the described top-down decomposition of the combinatorial objects makes these structures amenable to other methods, in particular, it seems that the contraction method, see, e.g., \cite{NeiRue2004,NeiSul2015}, allows an alternative characterization of limiting distributions occurring in the analysis of binary series-parallel networks.
Moreover, the combinatorial approach presented is flexible enough to allow also a study of series-parallel networks generated by modifications of the edge-duplication rules, in particular, one could treat also a Bernoulli model with a ``preferential edge-duplication rule'', or a $b$-ary saturation model by encoding the growth process via other recursive tree structures (edge-coloured plane increasing trees and bucket-recursive trees with bucket size $b \ge 2$, respectively); the authors plan to comment on that in a journal version of this work.

\section{Series-parallel networks and description via recursive tree structures\label{sec:Growth_models_recursive_trees}}

\subsection{Bernoulli model}

In the Bernoulli model in step $1$ one starts with a single edge labelled $1$ connecting the source and the sink, and in step $n$, with $n > 1$, one of the $n-1$ edges of the already generated series-parallel network is chosen uniformly at random, let us assume it is edge $j=(x,y)$; then either with probability $p$, $0 < p < 1$, this edge is doubled in a parallel way\footnote{In the original work \cite{Mahmoud2013} the r\^oles of $p$ and $q$ are switched, but we find it catchier to use $p$ for the probability of a parallel doubling.}, i.e., an additional edge $(x,y)$ labelled $n$ is inserted into the graph (let us say, right to edge $e$), or otherwise, thus with probability $q=1-p$, this edge is doubled in a serial way, i.e., edge $(x,y)$ is replaced by the series of edges $(x,z)$ and $(z,y)$, with $z$ a new node, where $(x,z)$ gets the label $j$ and $(z,y)$ will be labelled by $n$.

\begin{figure}
\begin{center}
\begin{minipage}{13cm}
\includegraphics[height=2.5cm]{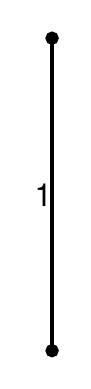}
\includegraphics[height=2.5cm]{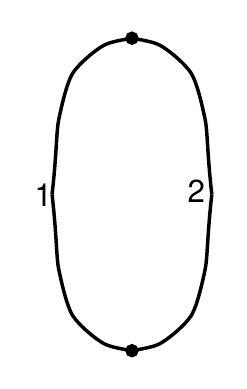}
\includegraphics[height=2.5cm]{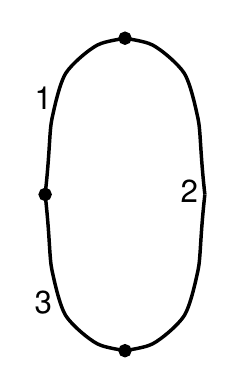}
\includegraphics[height=2.5cm]{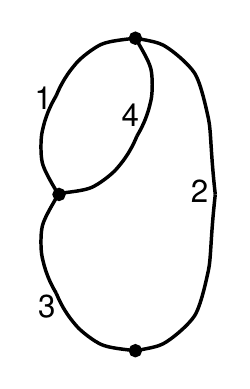}
\includegraphics[height=2.5cm]{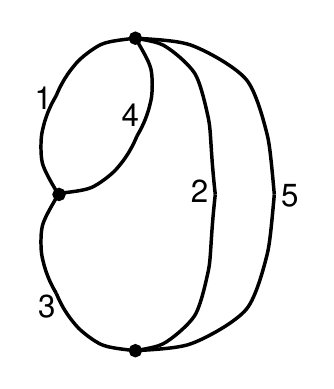}
\includegraphics[height=2.5cm]{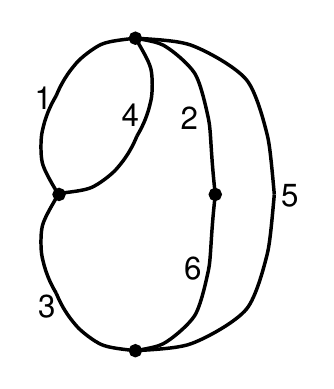}
\includegraphics[height=2.5cm]{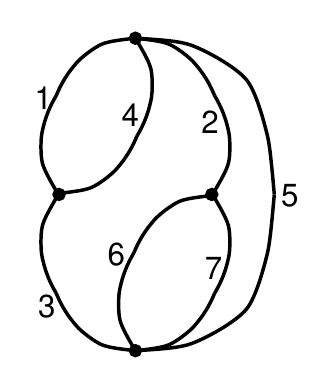}\hfill
\\
\includegraphics[height=2.5cm]{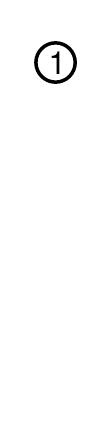}
\hspace*{5mm}\includegraphics[height=2.5cm]{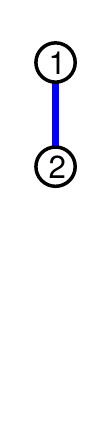}
\hspace*{5mm}\includegraphics[height=2.5cm]{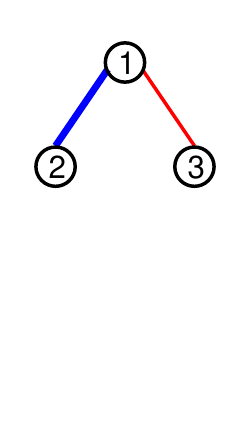}
\includegraphics[height=2.5cm]{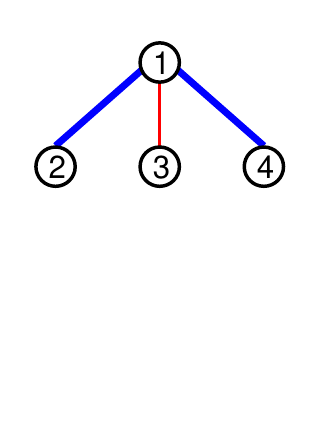}
\includegraphics[height=2.5cm]{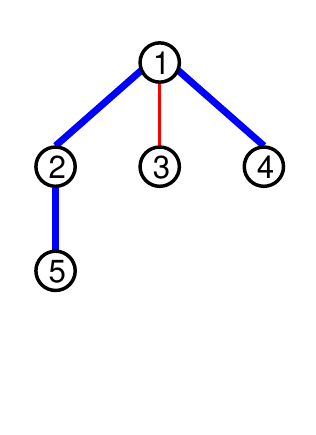}
\includegraphics[height=2.5cm]{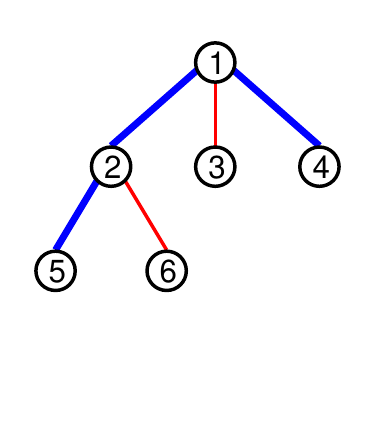}
\includegraphics[height=2.5cm]{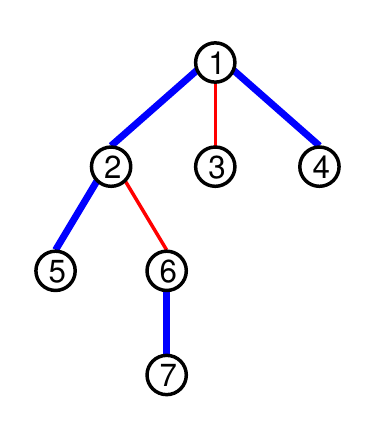}\hfill
\end{minipage}
\end{center}
\vspace*{-6mm}
\caption{Growth of a series-parallel network under the Bernoulli model and of the corresponding edge-coloured recursive tree. In the resulting graph the degree of the source is $4$, the length of the leftmost source-to-sink path is $2$ and there are $5$ different source-to-sink paths.\label{fig:growth_Bernoulli}}
\end{figure}

The growth of series-parallel networks corresponds to the growth of random recursive trees, where one starts in step $1$ with a node labelled $1$, and in step $n$ one of the $n-1$ nodes is chosen uniformly at random and node $n$ is attached to it as a new child. Thus, a doubling of edge $j$ in step $n$ when generating the series-parallel network corresponds in the recursive tree to an attachment of node $n$ to node $j$. Additionally, in order to keep the information about the kind of duplication of the chosen edge, the edge incident to $n$ is coloured either blue encoding a parallel doubling, or coloured red encoding a serial doubling. Such combinatorial objects of edge-coloured recursive trees can be described via the formal equation
\begin{equation*}
  \mathcal{T} = \mathcal{Z}^{\Box} \ast \text{\textsc{SET}}(B \cdot \mathcal{T} + R \cdot \mathcal{T}),
\end{equation*}
with $B$ and $R$ markers (see \cite{FlaSed2009}). Of course, one has to keep track of the number of blue and red edges to get the correct probability model according to 
\begin{equation*}
  \mathbb{P}\{\text{$T \in \mathcal{T}_{n}$ is chosen}\} = \frac{p^{\# \text{blue edges of $T$}} \cdot q^{\# \text{red edges of $T$}}}{T_{n}},
\end{equation*}
where $\mathcal{T}_{n} = \{T \in \mathcal{T} : \text{$T$ has order $n$}\}$ and $T_{n} := |\mathcal{T}_{n}| = (n-1)!$. Throughout this work the term order of a tree $T$ shall denote the number of labels contained in $T$, which, of course, for recursive trees coincides with the number of nodes of $T$. Then, each edge-coloured recursive tree of order $n$ and the corresponding series-parallel network of size $n$ occur with the same probability. An example for a series-parallel network grown via the Bernoulli model and the corresponding edge-coloured recursive tree is given in Figure~\ref{fig:growth_Bernoulli}.

\subsection{Binary model}

In the binary model again in step $1$ one starts with a single edge labelled $1$ connecting the source and the sink, and in step $n$, with $n>1$ one of the $n-1$ edges of the already generated series-parallel network is chosen uniformly at random; let us assume it is edge $j=(x,y)$; but now whether edge $j$ is doubled in a parallel or serial way is already determined by the out-degree of node $x$: if node $x$ has out-degree $1$ then we carry out a parallel doubling by inserting an additional edge $(x,y)$ labelled $n$ into the graph right to edge $j$, but otherwise, i.e., if node $x$ has out-degree $2$ and is thus already saturated, then we carry out a serial doubling by replacing edge $(x,y)$ by the edges $(x,z)$ and $(z,y)$, with $z$ a new node, where $(x,z)$ gets the label $j$ and $(z,y)$ will be labelled by $n$.

\begin{figure}
\begin{center}
\begin{minipage}{14cm}
\includegraphics[height=2.5cm]{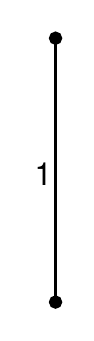}
\includegraphics[height=2.5cm]{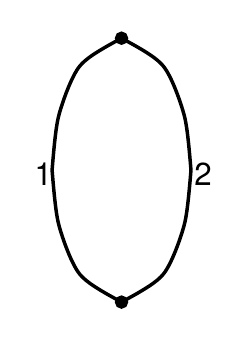}
\includegraphics[height=2.5cm]{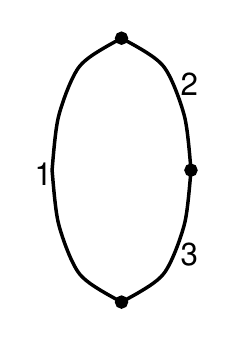}
\includegraphics[height=2.5cm]{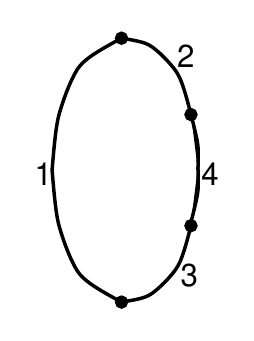}
\includegraphics[height=2.5cm]{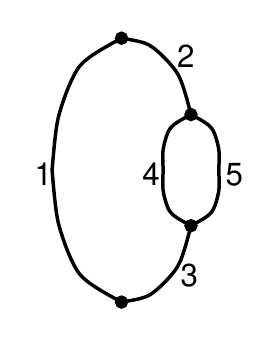}
\includegraphics[height=2.5cm]{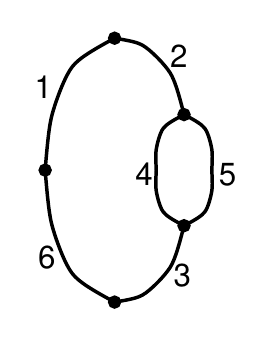}
\hspace*{4mm}\includegraphics[height=2.5cm]{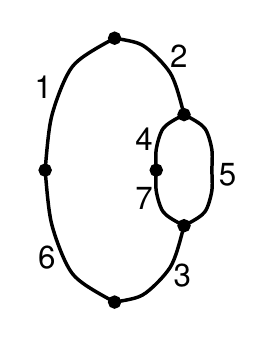}\hfill
\\
\includegraphics[height=2.5cm]{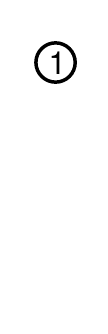}
\hspace*{3mm}\includegraphics[height=2.5cm]{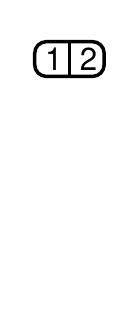}
\hspace*{6mm}\includegraphics[height=2.5cm]{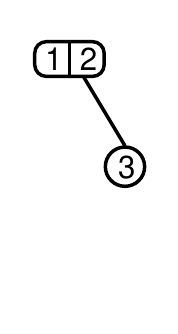}
\hspace*{3mm}\includegraphics[height=2.5cm]{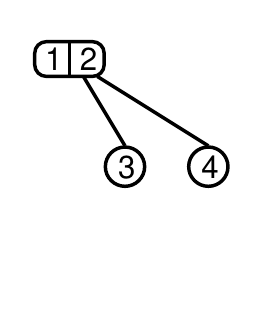}
\hspace*{-3mm}\includegraphics[height=2.5cm]{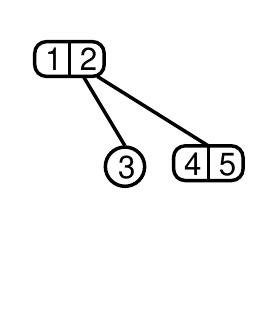}
\hspace*{-3mm}\includegraphics[height=2.5cm]{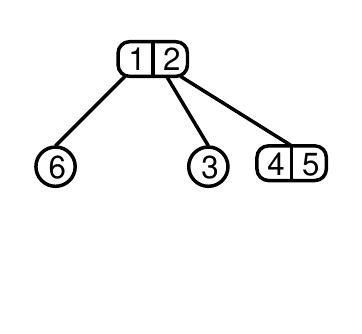}
\hspace*{-3mm}\includegraphics[height=2.5cm]{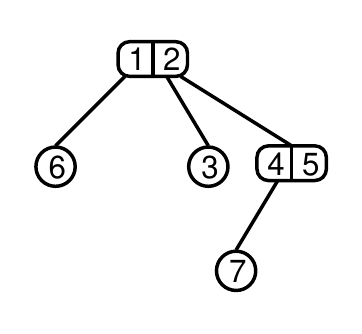}\hfill
\end{minipage}
\end{center}
\vspace*{-6mm}
\caption{Growth of a binary series-parallel network and of the corresponding bucket-recursive tree. In the resulting graph the degree of the sink is $2$, the length of the leftmost source-to-sink path is $2$ and there are $3$ different source-to-sink paths.\label{fig:growth_binary}}
\end{figure}

It turns out that the growth model for binary series-parallel networks corresponds with the growth model for bucket-recursive trees of maximal bucket size $2$, i.e., where nodes in the tree can hold up to two labels: in step $1$ one starts with the root node containing label $1$, and in step $n$ one of the $n-1$ labels in the tree is chosen uniformly at random, let us assume it is label $j$, and attracts the new label $n$. If the node $x$ containing label $j$ is saturated, i.e., it contains already two labels, then a new node containing label $n$ will be attached to $x$ as a new child, otherwise, label $n$ will be inserted into node $x$, then containing the labels $j$ and $n$. As has been pointed out in \cite{KubPan2010} such random bucket-recursive trees can also be described in a combinatorial way by extending the notion of increasing trees: namely a bucket-recursive tree is either a node labelled $1$ or it consists of the root node labelled $(1,2)$, where two (possibly empty) forests of (suitably relabelled) bucket-recursive trees are attached to the root as a left forest and a right forest. A formal description of the family $\mathcal{B}$ of bucket-recursive trees (of bucket size at most $2$) is in modern notation given as follows:
\begin{equation*}
  \mathcal{B} = \mathcal{Z}^{\Box} + \mathcal{Z}^{\Box} \ast \left(\mathcal{Z}^{\Box} \ast \left(\text{\textsc{SET}}(\mathcal{B}) \ast \text{\textsc{SET}}(\mathcal{B})\right)\right).
\end{equation*}
It follows from this formal description that there are $T_{n} = (n-1)!$ different bucket-recursive trees with $n$ labels, i.e., of order $n$, and furthermore it has been shown in \cite{KubPan2010} that this combinatorial description (assuming the uniform model, where each of these trees occurs with the same probability) indeed corresponds to the stochastic description of random bucket-recursive trees of order $n$ given before. An example for a binary series-parallel network and the corresponding bucket-recursive tree is given in Figure~\ref{fig:growth_binary}.

\begin{figure}
\begin{center}
\parbox{6.2cm}{\includegraphics[width=6cm]{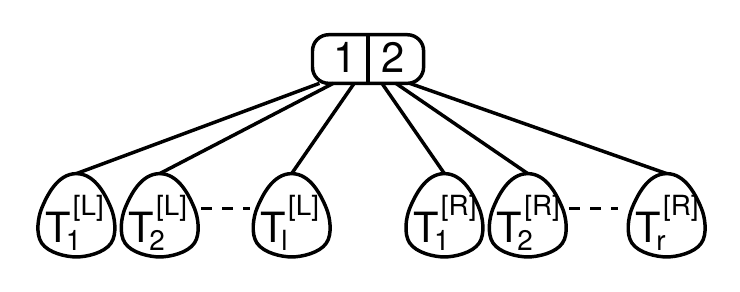}} \quad $\Longleftrightarrow$ \quad
\parbox{2.7cm}{\includegraphics[width=2.5cm]{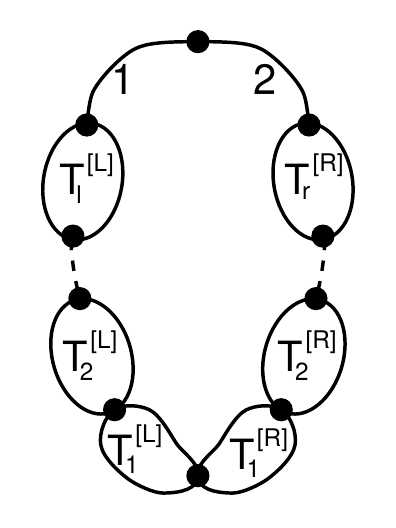}}
\end{center}
\vspace*{-3mm}
\caption{Decomposition of a bucket recursive tree $T$ into its root and the left and right forest, respectively, and the subblock-structure of the corresponding binary network.\label{fig:link_decomposition_buckettree_network}}
\end{figure}

In our analysis of binary series-parallel networks the following link between the decomposition of a bucket-recursive tree $T$ into its root $(1,2)$ and the left forest (consisting of the trees $T_{1}^{[L]}, \dots, T_{\ell}^{[L]}$) and the right forest (consisting of the trees $T_{1}^{[R]}, \dots, T_{r}^{[R]}$), and the subblock-structure of the corresponding binary network $G$ is important: $G$ consists of a left half $G^{[L]}$ and a right half $G^{[R]}$ (which share the source and the sink), where $G^{[L]}$ is formed by a series of blocks (i.e., maximal $2$-connected components) consisting of the edge labelled $1$ followed by binary networks corresponding to $T_{\ell}^{[L]}$, $T_{\ell-1}^{[L]}$, \dots, $T_{1}^{[L]}$, and $G^{[R]}$ is formed by a series of blocks consisting of the edge labelled $2$ followed by binary networks corresponding to $T_{r}^{[R]}$, $T_{r-1}^{[R]}$, \dots, $T_{1}^{[R]}$; see Figure~\ref{fig:link_decomposition_buckettree_network} for an example.

\section{Uniform Bernoulli edge-duplication growth model}

\subsection{Degree of the source\label{ssec:source-degree}}

Let $D_{n}=D_{n}(p)$ denote the r.v.\ measuring the degree of the source in a random series-parallel network of size $n$ for the Bernoulli model, with $0 < p < 1$. A first analysis of this quantity has been given in \cite{Mahmoud2013}, where the exact distribution of $D_{n}$ as well as exact and asymptotic results for the expectation $\mathbb{E}(D_{n})$ could be obtained. However, questions concerning the limiting behaviour of $D_{n}$ and the asymptotic behaviour of higher moments of $D_{n}$ have not been touched; in this context we remark that the explicit results for the probabilities $\mathbb{P}\{D_{n}=m\}$ as obtained in \cite{Mahmoud2013} and restated in \eqref{eqn:exact_dist_Dn} are not easily amenable for asymptotic studies, because of large cancellations of the alternating summands in the corresponding formula. We will reconsider this problem by applying the combinatorial approach introduced in Section~\ref{sec:Growth_models_recursive_trees}, and in order to get limiting distribution results we apply methods from analytic combinatorics. As has been already remarked in \cite{Mahmoud2013} the degree of the sink is equally distributed as $D_{n}$ due to symmetry reasons, although a simple justification of this fact via direct ``symmetry arguments'' does not seem to be completely trivial (the insertion process itself is a priori not symmetric w.r.t.\ the poles, since edges are always inserted towards the sink); however, it is not difficult to show this equality by establishing and treating a recurrence for the distribution of the sink, which is here omitted.

When considering the description of the growth process of these series-parallel networks via edge-coloured recursive trees it is apparent that the degree of the source in such a graph corresponds to the order of the maximal subtree containing the root node and only blue edges, i.e., we have to count the number of nodes in the recursive tree that can be reached from the root node by taking only blue edges; for simplicity we denote this maximal subtree by ``blue subtree''. Thus, in the recursive tree model, $D_{n}$ measures the order of the blue subtree in a random edge-coloured recursive tree of order $n$. To treat $D_{n}$ we introduce the r.v.\ $D_{n,k}$, whose distribution is given as the conditional distribution $D_{n} \big| \{\text{the tree has exactly $k$ blue edges}\}$, and the trivariate generating function
\begin{equation}
  F(z,u,v) := \sum_{n} \sum_{k} \sum_{m} T_{n} \binom{n-1}{k} \mathbb{P}\{D_{n,k} = m\} \frac{z^{n}}{n!} u^{k} v^{m},
\end{equation}
with $T_{n} = (n-1)!$ the number of recursive trees of order $n$. Thus $T_{n} \binom{n-1}{k} \mathbb{P}\{D_{n,k} = m\}$ counts the number of edge-coloured recursive trees of order $n$ with exactly $k$ blue edges, where the blue subtree has order $m$. Additionally we introduce the auxiliary function $N(z,u) := \sum_{n} \sum_{k} T_{n} \binom{n-1}{k} \frac{z^{n}}{n!} u^{k} = \frac{1}{1+u} \log\left(\frac{1}{1-z(1+u)}\right)$, i.e., the exponential generating function of the number of edge-coloured recursive trees of order $n$ with exactly $k$ blue edges.

The decomposition of a recursive tree into its root node and the set of branches attached to it immediately can be translated into a differential equation for $F(z,u,v)$, where we only have to take into account that the order of the blue subtree in the whole tree is one (due to the root node) plus the orders of the blue subtrees of the branches which are connected to the root node by a blue edge (i.e., only branches which are connected to the root node by a blue edge will contribute). Namely, with $F:= F(z,u,v)$ and $N := N(z,u)$, we get the first order separable differential equation
\begin{equation}\label{eqn:DEQ_source-degree}
  F' = v \cdot e^{u F + N},
\end{equation}
with initial condition $F(0,u,v) = 0$. Throughout this work, the notation $f'$ for (multivariate) functions $f(z, \dots)$ shall always denote the derivative w.r.t.\ the variable $z$. The exact solution of \eqref{eqn:DEQ_source-degree} can be obtained by standard means and is given as follows:
\begin{equation}
  F(z,u,v) = \frac{1}{u} \log \left(\frac{1}{1-v+v(1-z(1+u))^{\frac{u}{1+u}}}\right).
\end{equation}

Since we are only interested in the distribution of $D_{n}$ we will actually consider the generating function
\begin{equation}
  F(z,v) := \sum_{n} \sum_{m} T_{n} \mathbb{P}\{D_{n} = m\} \frac{z^{n}}{n!} v^{m} = \sum_{n} \sum_{m} \mathbb{P}\{D_{n} = m\} \frac{z^{n}}{n} v^{m}.
\end{equation}
According to the definition of the conditional r.v.\ $D_{n,k}$ it holds that $\mathbb{P}\{D_{n}=m\} = \sum_{k=0}^{n-1} \mathbb{P}\{D_{n,k}=m\} \binom{n-1}{k} p^{k} q^{n-1-k}$, which, after simple computations, gives the relation $F(z,v) = \frac{1}{q} F(qz, \frac{p}{q}, v)$. Thus we obtain the following explicit formula for $F'(z,v)$, which has been obtained already in \cite{Mahmoud2013} by using a description of $D_{n}$ via urn models:
\begin{equation}\label{eqn:gf_sourcedeg_explicit}
  F'(z,v) = \frac{v}{v(1-z) + (1-v) (1-z)^{1-p}}.
\end{equation}
Extracting coefficients from \eqref{eqn:gf_sourcedeg_explicit} immediately yields the explicit result for the probability distribution of $D_{n}$, with $1 \le m \le n$, stated in \cite{Mahmoud2013}:
\begin{equation}\label{eqn:exact_dist_Dn}
  \mathbb{P}\{D_{n}=m\} = [z^{n-1} v^{m}] F'(z,v) = \sum_{j=0}^{m-1} \binom{m-1}{j} (-1)^{n+j-1} \binom{p(j+1)-1}{n-1}.
\end{equation}

In order to describe the limiting distribution behaviour of $D_{n}$ we first study the integer moments. To do this we introduce $\tilde{F}(z,w) := F(z,1+w)$, since we get for its derivative the relation $\tilde{F}'(z,w) = \sum_{n} \sum_{r} \mathbb{E}(D_{n}^{\underline{r}}) z^{n-1} \frac{w^{r}}{r!}$, with $\mathbb{E}(D_{n}^{\underline{r}}) = \mathbb{E}(D_{n} \cdot (D_{n}-1) \cdots (D_{n}-r+1))$ the $r$-th factorial moment of $D_{n}$. Plugging $v=1+w$ into \eqref{eqn:gf_sourcedeg_explicit}, extracting coefficients and applying Stirling's formula for the factorials easily gives the following explicit and asymptotic result for the $r$-th factorial moments of $D_{n}$, with $r \ge 1$:
\begin{equation*}
  \mathbb{E}(D_{n}^{\underline{r}}) = r! \sum_{j=0}^{r-1} \binom{r-1}{j} (-1)^{r-1-j} \binom{n+p(j+1)-1}{n-1} \sim 
	\frac{r! \cdot n^{rp}}{\Gamma(rp+1)},
\end{equation*}
from which we further deduce
\begin{equation}
  \mathbb{E}\left(\big(\frac{D_{n}}{n^{p}}\big)^{r}\right) \sim \frac{r!}{\Gamma(rp+1)}.
\end{equation}
Thus, the $r$-th integer moments of the suitably scaled r.v.\ $D_{n}$ converge to the integer moments of a so-called Mittag-Leffler distribution $D=D(p)$ with parameter $p$ (see, e.g., \cite{Janson2010}), which, by an application of the theorem of Fr\'echet and Shohat, indeed characterizes the limiting distribution of $D_{n}$.

From the explicit formula \eqref{eqn:gf_sourcedeg_explicit} it is also possible to characterize the density function $f(x)$ of $D$ (We remark that alternatively we can obtain $f(x)$ from the moment generating function $M(z) = \mathbb{E}(e^{D z}) = \sum_{r \ge 0} \mathbb{E}(D^{r}) \frac{z^{r}}{r!}$ and applying the inverse Laplace transform.). Namely, it holds
\begin{equation}\label{eqn:probability_Hankel}
  \mathbb{P}\{D_{n}=m\} = [z^{n-1}v^{m}] F'(z,v) = \frac{1}{2 \pi i} \oint \frac{(1-(1-z)^{p})^{m-1}}{z^{n} (1-z)^{1-p}} dz,
\end{equation}
where we have to choose as contour a positively oriented simple closed curve around the origin, which lies in the domain of analyticity of the integrand. To evaluate the integral asymptotically (and uniformly) for $m = O(n^{p+\delta})$, $\delta > 0$ and $n \to \infty$ one can adapt the considerations done in \cite{Panholzer2004} for the particular instance $p=\frac{1}{2}$. After straightforward computations one obtains the following asymptotic equivalent of these probabilities, which determines the density function $f(x)$ of the limiting distribution $D$:
\begin{equation*}
  \mathbb{P}\{D_{n} = m\} \sim \frac{1}{n^{p}} \cdot \frac{1}{2\pi i} \int_{\mathcal{H}} \frac{e^{-t-\frac{m}{n^{p}} (-t)^{p}}}{(-t)^{1-p}} dt,
\end{equation*}
with $\mathcal{H}$ a Hankel contour starting from $e^{2 \pi i} \infty$, passing around $0$ and terminating at $+ \infty$.
The results are collected in the following theorem.
\begin{theorem}\label{thm:Bernoulli_degree_limit}
The degree $D_{n}$ of the source or the sink in a randomly chosen series-parallel network of size $n$ generated by the Bernoulli model converges after scaling, for $n \to \infty$, in distribution to a Mittag-Leffler distribution $D=D(p)$ with parameter $p$:
$\frac{D_{n}}{n^{p}} \xrightarrow{(d)} D$, where $D$ is characterized by the sequence of its $r$-th integer moments:
\begin{equation*}
  \mathbb{E}(D^{r}) = \frac{r!}{\Gamma(rp+1)}, \quad \text{for $r \ge 0$},
\end{equation*}
as well as by its density function $f(x)$ (with $\mathcal{H}$ a Hankel contour):
\begin{equation*}
  f(x) = \frac{1}{2 \pi i} \int_{\mathcal{H}} \frac{e^{-t - x (-t)^{p}}}{(-t)^{1-p}} dt, \quad \text{for $x > 0$}.
\end{equation*}
\end{theorem}

We remark that after simple manipulations we can also write $f(x)$ as the following real integral:
\begin{equation*}
  f(x) = \frac{1}{\pi p} \int_{0}^{\infty} e^{-w^{\frac{1}{p}} - xw \cos(\pi p)} \sin(\pi p - x w \sin(\pi p)) dw, \quad \text{for $x>0$}.
\end{equation*}
We further remark that for the particular instance $p=\frac{1}{2}$ one can evaluate the Hankel integral above and obtains that the limiting distribution $D$ is characterized by the density function $f(x) = \frac{1}{\sqrt{\pi}} \cdot e^{-\frac{x^{2}}{4}}$, $x > 0$. Thus, $f(x)$ is the density function of a so-called half-normal distribution with parameter $\sigma = \sqrt{2}$.

\subsection{Length of a random path from source to sink}

We consider the length $L_{n} = L_{n}(p)$ (measured by the number of edges) of a random path from the source to the sink in a randomly chosen series-parallel network of size $n$ for the Bernoulli model. In this context, the following definition of a random source-to-sink path seems natural: we start at the source and walk along outgoing edges, such that whenever we reach a node of out-degree $d$, $d \ge 1$, we choose one of these outgoing edges uniformly at random, until we arrive at the sink.

The following two observations are very helpful in the analysis of this parameter. First it holds that the length $L_{n}$ of a random path is distributed as the length $L_{n}^{[L]}$ of the leftmost source-to-sink path in a random series-parallel network of size $n$; the meaning of the leftmost path is, that whenever we reach a node of out-degree $d$, we choose the first (i.e., leftmost) outgoing edge. Unfortunately, so far we do not see a simple symmetry argument to show this fact (such an argument easily shows that the rightmost path has the same distribution as the leftmost path, but it does not seem to explain the general situation). However, we are able to show this in a somehow indirect manner: namely, it is possible to establish a more involved distributional recurrence for the length of a random path $L_{n}$ and show that the explicit solution for the probability distribution of the length of the leftmost path $L_{n}^{[L]}$ is indeed the solution of the recurrence for $L_{n}$. These computations will be given in the journal version of this work, here we have to omit them.

Second we use that the length of the leftmost source-to-sink path in a series-parallel network has the following simple description in the corresponding edge-coloured recursive tree: namely, an edge is lying on the leftmost source-to-sink path if and only if the corresponding node in the recursive tree can be reached from the root by using only red edges (i.e., edges that correspond to serial edges). This means that the length $\ell$ of the leftmost source-to-sink path corresponds in the edge-coloured recursive tree model to the order of the maximal subtree containing the root node and only red edges. If we switch the colours red and blue in the tree we obtain an edge-coloured recursive tree where the maximal blue subtree has the same order, i.e., where the source-degree of the corresponding series-parallel network is $\ell$. But switching colours in the tree model corresponds to switching the probabilities $p$ and $q=1-p$ for generating a parallel and a serial edge, respectively, in the series-parallel network. Thus it simply holds $L_{n}^{[L]}(p) \stackrel{(d)}{=} D_{n}(1-p)$, where $D_{n}$ denotes the source-degree in a random series-parallel network of size $n$.

Combining these considerations and the results for $D_{n}$ obtained in Section~\ref{ssec:source-degree} we obtain the following theorem.
\begin{theorem}\label{thm:Bernoulli_length_path}
The length $L_{n}$ of a random path from the source to the sink in a randomly chosen series-parallel network of size $n$ generated by the Bernoulli model has the following probability distribution:
\begin{equation*}
  \mathbb{P}\{L_{n}=m\} = \sum_{j=0}^{m-1} \binom{m-1}{j} (-1)^{n+j-1} \binom{j-p(j+1)}{n-1}, \quad \text{for $1 \le m \le n$}.
\end{equation*}
Moreover, $L_{n}$ has, for $n \to \infty$, the following limiting distribution behaviour:
$\frac{L_{n}}{n^{1-p}} \xrightarrow{(d)} L$, where the limiting distribution $L$ is a Mittag-Leffler distribution with parameter $1-p$, i.e., $L$ is characterized by the sequence of its $r$-th integer moments:
\begin{equation*}
  \mathbb{E}(L^{r}) = \frac{r!}{\Gamma(r(1-p)+1)}, \quad \text{for $r \ge 0$},
\end{equation*}
as well as by its density function $g(x)$ (with $\mathcal{H}$ a Hankel contour):
\begin{equation*}
  g(x) = \frac{1}{2 \pi i} \int_{\mathcal{H}} \frac{e^{-t - x (-t)^{1-p}}}{(-t)^{p}} dt, \quad \text{for $x > 0$},
\end{equation*}
\end{theorem}

\subsection{Number of paths from source to sink}

Let $P_{n} = P_{n}(p)$ denote the r.v.\ measuring the number of different paths from the source to the sink in a randomly chosen series-parallel network of size $n$ for the Bernoulli model. Again we use the description of the growth of the graphs via edge-labelled recursive trees, but in contrast to the previous studies of parameters, here it seems advantageous to use an alternative decomposition of recursive trees with respect to the edge connecting nodes $1$ and $2$, which allows to establish a stochastic recurrence for the r.v.\ $P_{n}$. Namely, it is not difficult to show (see, e.g., \cite{DobFil1999}) that when starting with a random recursive tree $T$ of order $n \ge 2$ and removing the edge $1-2$, both resulting trees $T'$ and $T''$ are (after an order-preserving relabelling) again random recursive trees of smaller orders; moreover, if $U_{n}$ denotes the order of the resulting tree $T'$ rooted at the former label $2$ (and thus $n-U_{n}$ gives the order of the tree $T''$ rooted at the original root of the tree $T$), it holds that $U_{n}$ follows a discrete uniform distribution on the integers $\{1, \dots, n-1\}$, i.e., $\mathbb{P}\{U_{n} = k\} = \frac{1}{n-1}$, for $1 \le k \le n-1$. Depending on the colour of the edge $1-2$ in the edge-labelled recursive tree, it corresponds to a parallel edge (colour blue, which occurs with probability $p$) or a serial edge (colour red, which occurs with probability $q=1-p$) in the series-parallel network: if it is a parallel edge then the number of source-to-sink paths in the corresponding substructures have to be added, whereas for a serial edge they have to be multiplied in order to obtain the total number of source-to-sink paths in the whole graph. Thus $P_{n}$ satisfies the following stochastic recurrence:
\begin{equation}\label{eqn:nrpath_stochrec}
  P_{n} \stackrel{(d)}{=} \boldsymbol{1}_{\{B_{n}=1\}} \cdot \left(P_{U_{n}}' + P_{n-U_{n}}''\right) + \boldsymbol{1}_{\{B_{n}=0\}} \cdot \left(P_{U_{n}}' \cdot P_{n-U_{n}}''\right), \quad \text{for $n \ge 2$}, \quad P_{1} = 1,
\end{equation}
where $B_{n}$ and $U_{n}$ are independent of each other and independent of $(P_{n})$, $(P_{n}')$ and $(P_{n}'')$, and where $(P_{n}')$ and $(P_{n}'')$ are independent copies of $(P_{n})$. Here $B_{n}$ is the indicator variable of the event that $1-2$ is a blue edge in the recursive tree, thus $B_{n}$ is a Bernoulli distributed r.v.\ with success probability $p$, i.e., $\mathbb{P}\{B_{n} = 1\} = p$. Furthermore, the r.v.\ $U_{n}$ measuring the order of the subtree rooted at $2$, is uniformly distributed on $\{1, 2, \dots, n-1\}$, i.e., $\mathbb{P}\{U_{n} = k\} = \frac{1}{n-1}$, for $1 \le k \le n-1$.

Starting with \eqref{eqn:nrpath_stochrec} and taking the expectations yields after simple manipulations the following recurrence:
\begin{equation}
\begin{split}
  \mathbb{E}(P_{n}) & = \frac{2p}{n-1} \sum_{k=1}^{n-1} \mathbb{E}(P_{k}) + \frac{1-p}{n-1} \sum_{k=1}^{n-1} \mathbb{E}(P_{k}) \mathbb{E}(P_{n-k}), \quad n \ge 2, \qquad \mathbb{E}(P_{1})=1.
\end{split}
\end{equation}
To treat this recurrence we introduce the generating function $E(z) := \sum_{n \ge 1} \mathbb{E}(P_{n}) z^{n-1}$, which gives the following non-linear first order differential equation of Bernoulli type:
\begin{equation}\label{eqn:DEQ_Bernoulli_nrpath}
  E'(z) = \frac{2p}{1-z} E(z) + (1-p) \big(E(z)\big)^{2}, \quad E(0)=1.
\end{equation}

Equation~\eqref{eqn:DEQ_Bernoulli_nrpath} can be treated by a standard technique for Bernoulli type differential equations and leads to the following solution, where we have to distinguish whether $p=\frac{1}{2}$ or not:
\begin{equation}\label{eqn:GF_exp_Ez}
  E(z) = 
  \begin{cases}
    \frac{1-2p}{(1-p)(1-z) - p(1-z)^{2p}}, & \quad \text{for $p \neq \frac{1}{2}$},\\
    \frac{2}{2(1-z) - (1-z)\log\left(\frac{1}{1-z}\right)}, & \quad \text{for $p = \frac{1}{2}$}.
	\end{cases}
\end{equation}

From the formula \eqref{eqn:GF_exp_Ez} for the generating function $E(z)$ one can easily deduce explicit results for the expected value $\mathbb{E}(P_{n}) = [z^{n-1}] E(z)$, which, however, due to alternating signs of the summands are not easily amenable for asymptotic considerations. Instead, in order to obtain the asymptotic behaviour of $\mathbb{E}(P_{n})$ we consider the formul{\ae} for the generating function $E(z)$ stated in \eqref{eqn:GF_exp_Ez} and describe the structure of the singularities: for $0 < p < 1$ the dominant singularity at $z=\rho<1$ is annihilating the denominator; there $E(z)$ has a simple pole, which due to singularity analysis \cite{FlaSed2009} yields the main term of $\mathbb{E}(P_{n})$, i.e., the asymptotically exponential growth behaviour; the (algebraic or logarithmic) singularity at $z=1$ determines the second and higher order terms in the asymptotic behaviour of $\mathbb{E}(P_{n})$, which differ according to the ranges $0<p<\frac{1}{2}$, $p=\frac{1}{2}$, and $\frac{1}{2} < p < 1$.
This yields the following theorem.
\begin{theorem}
The expectation $\mathbb{E}(P_{n})$ of the number of paths $P_{n}$ from source to sink in a random series-parallel network of size $n$ generated by the Bernoulli model is given by the following explicit formula:
\begin{equation*}
  \mathbb{E}(P_{n}) = 
	\begin{cases}
	  \sum_{j=0}^{n-1} (-1)^{n+j-1} \binom{(2p-1)j-1}{n-1} \sum_{k=0}^{n-1} \binom{k}{j} \left(\frac{p}{2p-1}\right)^{k}, & \quad \text{for $p \neq \frac{1}{2}$},\\
		\sum_{k=0}^{n-1} \frac{(-1)^{k}}{2^{k}} \cdot B_{k}(-H_{n-1}^{(1)},-H_{n-1}^{(2)},-2H_{n-1}^{(3)},\ldots,-(k-1)! H_{n-1}^{(k)}), & \quad \text{for $p=\frac{1}{2}$},
	\end{cases}
\end{equation*}
where $B_{k}(x_{1}, x_{2}, \ldots, x_{k})$ denotes the $k$-th complete Bell polynomial and where $H_{n}^{(m)} := \sum_{j=1}^{n} \frac{1}{j^{m}}$ denote the $m$-th order harmonic numbers (see \cite{Zave1976}).

The asymptotic behaviour of $\mathbb{E}(P_{n})$ is, for $n \to \infty$, given as follows:
\begin{gather*}
  \mathbb{E}(P_{n}) = \frac{1}{1-p} \cdot \alpha_{p}^{n} + R_{p}(n),\\
\begin{split}
& \text{where $\alpha_{p} = \frac{1}{1-\big(\frac{p}{1-p}\big)^{\frac{1}{1-2p}}}$, for $p \neq \frac{1}{2}$, \enspace and \enspace
$\alpha_{p} = \frac{1}{1-e^{-2}} = \lim_{p \to \frac{1}{2}}\frac{1}{1-\big(\frac{p}{1-p}\big)^{\frac{1}{1-2p}}}$, for $p = \frac{1}{2}$},\\
& \text{and $R_{p}(n) = - \frac{1-2p}{p \Gamma(2p)} n^{2p-1} + \mathcal{O}(n^{2(2p-1)})$, for $0<p<\frac{1}{2}$, \enspace 
$R_{p}(n) = -\frac{2}{\log n} + \mathcal{O}(\frac{1}{\log^{2}n})$, for $p=\frac{1}{2}$},\\ 
& \quad \text{$R_{p}(n) = - \frac{2p-1}{1-p} + \mathcal{O}(n^{1-2p})$, for $\frac{1}{2}<p<1$}.
\end{split}
\end{gather*}
\end{theorem}

\section{Uniform binary saturation edge-duplication growth model}

\subsection{Length of a random path from source to sink}

We are interested in the length of a typical source-to-sink path in a series-parallel network of size $n$. Again, it is natural to start at the source of the graph and move along outgoing edges, in a way that whenever we have the choice of two outgoing edges we use one of them uniformly at random to enter a new node, until we finally end at the sink. Let us denote by $L_{n}$ the length of such a random source-to-sink path in a random series-parallel network of size $n$ for the binary model.
Due to symmetry reasons it holds that $L_{n} \stackrel{(d)}{=} L_{n}^{[L]}$, where $L_{n}^{[L]}$ denotes the length of the leftmost source-to-sink path in a random series-parallel network of size $n$, i.e., the source-to-sink path, where in each node we choose the left outgoing edge to enter the next node.

In order to analyse $L_{n}^{[L]}$ we use the description of the growth of series-parallel networks via bucket-recursive trees: the length of the left path is equal to $1$ (coming from the root node of the tree, i.e., stemming from the edge $1$ in the graph) plus the sum of the lengths of the left paths in the subtrees contained in the left forest (which correspond to the blocks of the left part of the graph). When we introduce the generating function
\begin{equation}
  F(z,v) := \sum_{n \ge 1} \sum_{m \ge 0} T_{n} \mathbb{P}\{L_{n}=m\} \frac{z^{n}}{n!} v^{m}
	= \sum_{n \ge 1} \sum_{m \ge 0} \mathbb{P}\{L_{n}=m\} \frac{z^{n}}{n} v^{m},
\end{equation}
then the above description yields the following differential equation:
\begin{equation}\label{eqn_DEQ_binary_length}
  F''(z,v) = v e^{F(z,v)} e^{N(z)} = \frac{v}{1-z} e^{F(z,v)}, \quad F(0,v)=0, \quad F'(0,v)=v,
\end{equation}
where $N(z) = \log\frac{1}{1-z}$ is the exponential generating function of the number $T_{n} = (n-1)!$ of bucket-recursive trees of order $n$. In order to compute the expectation we consider $E(z) := \left.\frac{\partial}{\partial v} F(z,v)\right|_{v=1} = \sum_{n \ge 1} \mathbb{E}(L_{n}) \frac{z^{n}}{n}$, which satisfies the following linear second order differential equation of Eulerian type:
\begin{equation*}
  E''(z) = \frac{1}{(1-z)^{2}} E(z) + \frac{1}{(1-z)^{2}}, \quad E(0)=0, \quad E'(0)=1.
\end{equation*}
The explicit solution can be obtained by a standard technique and is given as follows:
\begin{equation*}
  E(z) = \frac{3+\sqrt{5}}{2 \sqrt{5}} \frac{1}{(1-z)^{\frac{\sqrt{5}-1}{2}}} - \frac{3-\sqrt{5}}{2 \sqrt{5}} (1-z)^{\frac{1+\sqrt{5}}{2}} - 1.
\end{equation*}
Extracting coefficients and applying Stirling's formula immediately yields the following explicit and asymptotic result for the expectation:
\begin{equation}
  \mathbb{E}(L_{n}) = n \left(\frac{3+\sqrt{5}}{2 \sqrt{5}} \binom{n+\frac{\sqrt{5}}{2}-\frac{3}{2}}{n} - \frac{3-\sqrt{5}}{2 \sqrt{5}} \binom{n-\frac{\sqrt{5}}{2} - \frac{3}{2}}{n}\right) \sim \frac{1+\sqrt{5}}{2 \sqrt{5}} \frac{n^{\frac{\sqrt{5}-1}{2}}}{\Gamma(\frac{\sqrt{5}-1}{2})}.
\end{equation}

In order to characterize the limiting distribution of $L_{n}$ we will compute iteratively the asymptotic behaviour of all its integer moments. To this aim it is advantageous to consider $G(z,v) := F'(z,v)$. Differentiating \eqref{eqn_DEQ_binary_length} shows that $G(z,v)$ satisfies the following differential equation:
\begin{equation}\label{eqn:DEQ_binary_length_G}
  G''(z,v) = G'(z,v) G(z,v) + \frac{1}{1-z} G'(z,v), \quad G(0,v)=v, \quad G'(0,v) = v.
\end{equation}
Introducing $M_{r}(z) := \left.\frac{\partial^{r}}{\partial v^{r}} G(z,v)\right|_{v=1} = \sum_{n \ge 1} \mathbb{E}(L_{n}^{\underline{r}}) z^{n-1}$, differentiating \eqref{eqn:DEQ_binary_length_G} $r$ times w.r.t.\ $v$ and evaluating at $v=1$ yields
\begin{equation*}
  M_{r}''(z) = \frac{2}{1-z} M_{r}'(z) + \frac{1}{(1-z)^{2}} M_{r}(z) + R_{r}(z),
\end{equation*}
with $R_{r}(z) = \sum_{k=1}^{r-1} \binom{r}{k} M_{k}'(z) M_{r-k}(z)$.
Thus $M_{r}(z)$ satisfies for each $r$ an Eulerian differential equation, where the inhomogeneous part $R_{r}(z)$ depends on the functions $M_{k}(z)$, with $k < r$. The asymptotic behaviour of $M_{r}(z)$ around the dominant singularity $z=1$ can be established inductively, namely it holds:
\begin{equation*}
  M_{r}(z) \sim \frac{c_{r}}{(1-z)^{r \tilde{\phi}+1}},
\end{equation*}
with $\tilde{\phi} = \frac{\sqrt{5}-1}{2}$, and where the constants $c_{r}$ satisfy a certain recurrence of ``convolution type''.
Singularity analysis and an application of the theorem of Fr\'{e}chet and Shohat shows then the following limiting distribution result.
\begin{theorem}
  The length $L_{n}$ of a random path from the source to the sink in a random series-parallel network of size $n$ generated by the binary model satisfies, for $n \to \infty$, the following limiting distribution behaviour, with $\tilde{\phi}=\frac{\sqrt{5}-1}{2}$:
	\begin{equation*}
	  \frac{L_{n}}{n^{\tilde{\phi}}} \xrightarrow{(d)} L,
	\end{equation*}
	where the limiting distribution $L$ is characterized by its sequence of $r$-th integer moments via
	\begin{equation*}
	  \mathbb{E}(L^{r}) = \frac{r! \cdot \tilde{c}_{r}}{\Gamma(r \tilde{\phi} +1)}, \quad r \ge 0,
	\end{equation*}
	and where the sequence $\tilde{c}_{r}$ satisfies the recurrence
	$\tilde{c}_{r} = \frac{1}{\tilde{\phi} (r-1) ((r+1)\tilde{\phi}+1)} \sum_{k=1}^{r-1} (k\tilde{\phi}+1) \tilde{c}_{k} \tilde{c}_{r-k}$, 
	for $r \ge 2$, with $\tilde{c}_{0}=1$ and $\tilde{c}_{1} = \frac{3+\tilde{\phi}}{5}$.
\end{theorem}

\subsection{Degree of the sink}

Whereas the (out-)degree of the source of a binary series-parallel network is two (if the graph has at least two edges), typically the (in-)degree of the sink is quite large, as will follow from our treatments. Let us denote by $D_{n}$ the degree of the sink in a random series-parallel network of size $n$ for the binary model. For a binary series-parallel network, the value of this parameter can be determined recursively by adding the degrees of the sinks in the last block of each half of the graph; in the case that a half only consists of one edge then the contribution of this half is of course $1$. When considering the corresponding bucket-recursive tree this means that the degree of the sink can be computed recursively by adding the contributions of the left and the right forest attached to the root, where the contribution of a forest is either given by $1$ in case that the forest is empty (then the corresponding root node contributes to the degree of the sink) or it is the contribution of the first tree in the forest (which corresponds to the last block), see Figure~\ref{fig:link_decomposition_buckettree_network}.
Introducing the generating functions
\begin{equation}
  F(z,v) := \sum_{n \ge 1} \sum_{m \ge 0} T_{n} \mathbb{P}\{D_{n}=m\} \frac{z^{n}}{n!} v^{m}, \quad
	A(z,v) := \sum_{n \ge 0} \sum_{m \ge 0} \tilde{T}_{n} \mathbb{P}\{\tilde{D}_{n}=m\} \frac{z^{n}}{n!} v^{m},
\end{equation}
with $\tilde{D}_{n}$ denoting the corresponding quantity for the left or right forest and $\tilde{T}_{n} = n!$ counting the number of forests of order $n$,
the combinatorial decomposition of bucket-recursive trees yields the following system of differential equations:
\begin{equation}\label{eqn:Fzv_Azv}
  F''(z,v) = \big(A(z,v)\big)^{2}, \quad A'(z,v) = \frac{1}{1-z} \cdot F'(z,v).
\end{equation}
From system \eqref{eqn:Fzv_Azv} the following non-linear differential equation for $F(z,v)$ can be obtained:
\begin{equation*}
  F'''(z,v) = \frac{2}{1-z} \sqrt{F''(z,v)} F'(z,v), \quad F(0,v)=0, F'(0,v)v, F''(0,v)=v^{2},
\end{equation*}
which, by considering $E(z) := \left.\frac{\partial}{\partial v}F(z,v)\right|_{v=1}$ and solving an Eulerian differential equation, allows to compute an exact and asymptotic expression for the expectation; namely it holds
\begin{equation}
  \mathbb{E}(D_{n}) = \frac{1+\sqrt{2}}{2} \binom{n+\sqrt{2}-2}{n-1} - \frac{\sqrt{2}-1}{2} \binom{n-\sqrt{2}-2}{n-1} \sim \frac{1+\sqrt{2}}{2} \frac{n^{\sqrt{2}-1}}{\Gamma(\sqrt{2})}.
\end{equation}

However, for asymptotic studies of higher moments it seems to be advantageous to consider the following second order non-linear differential equation for $A(z,v)$, which follows immediately from \eqref{eqn:Fzv_Azv}:
\begin{equation}\label{eqn:DEQ_degree_binary_A}
  A''(z,v) = \frac{1}{1-z} A'(z,v) + \frac{1}{1-z} \big(A(z,v)\big)^{2}, \quad A(0,v)=v, \quad A'(0,v)=v.
\end{equation}
Introducing the functions $\tilde{M}_{r}(z) := \left.\frac{\partial^{r}}{\partial v^{r}}A(z,v)\right|_{v=1}$ and differentiating \eqref{eqn:DEQ_degree_binary_A} $r$ times, one obtains that $\tilde{M}_{r}(z)$ satisfies for $r \ge 1$ the following second order Eulerian differential equation:
\begin{equation}\label{eqn:DEQ_degree_binary_moments}
  \tilde{M}_{r}''(z) = \frac{1}{1-z} \tilde{M}_{r}'(z) + \frac{2}{(2-z)^{2}} \tilde{M}_{r}(z) + R_{r}(z),
\end{equation}
with $R_{r}(z) := \frac{1}{1-z} \sum_{k=1}^{r-1} \binom{r}{k} \tilde{M}_{k}(z) \tilde{M}_{r-k}(z)$.
From \eqref{eqn:DEQ_degree_binary_moments} one can inductively show that the local behaviour of the functions $\tilde{M}_{r}(z)$ around the (unique) dominant singularity $z=1$ is given as follows:
\begin{equation*}
  \tilde{M}_{r}(z) \sim \frac{c_{r}}{(1-z)^{r \sqrt{2} - (r-1)}}, \quad r \ge 1,
\end{equation*}
where the constants $c_{r}$ are determined recursively.
Actually we are interested in the functions $M_{r}(z) := \left.\frac{\partial^{r}}{\partial v^{r}}F(z,v)\right|_{v=1} = \sum_{n \ge 1} \mathbb{E}(D_{n}^{\underline{r}}) \frac{z^{n}}{n}$, which are, due to \eqref{eqn:Fzv_Azv}, related to $\tilde{M}_{r}(z)$ via $M_{r}''(z) = (1-z) \tilde{M}_{r}''(z) - \tilde{M}_{r}'(z)$. Singularity analysis as well as the theorem of Fr\'{e}chet and Shohat show then the following limiting distribution result.
\begin{theorem}
  The degree $D_{n}$ of the sink in a random series-parallel network of size $n$ generated by the binary model satisfies, for $n \to \infty$, the following limiting distribution behaviour:
	\begin{equation*}
	  \frac{D_{n}}{n^{\sqrt{2}-1}} \xrightarrow{(d)} D,
	\end{equation*}
	where the limiting distribution $D$ is characterized by its sequence of $r$-th integer moments via
	\begin{equation*}
	  \mathbb{E}(D^{r}) = \frac{r! (r(\sqrt{2}-1)+1)\tilde{c}_{r}}{\Gamma(r (\sqrt{2}-1) +1)}, \quad r \ge 0,
	\end{equation*}
	where the sequence $\tilde{c}_{r}$ satisfies the recurrence
	$\tilde{c}_{r} = \frac{1}{(r(\sqrt{2}-1)+1)^{2}-2} \sum_{k=1}^{r-1} \tilde{c}_{k} \tilde{c}_{r-k}$, for $r \ge 2$, with $\tilde{c}_{0}=1$ and $\tilde{c}_{1} = \frac{1+\sqrt{2}}{2 \sqrt{2}}$.
\end{theorem}

\subsection{Number of paths from source to sink}

As for the Bernoulli model we are interested in results concerning the number of different paths from the source to the sink in a series-parallel network and denote by $P_{n}$ the number of source-to-sink paths in a random series-parallel network of size $n$ for the binary model. In order to study $P_{n}$ it seems advantageous to start with a stochastic recurrence for this random variable obtained by decomposing the bucket-recursive tree into the root node and the left and right forest (of bucket-recursive trees) attached to the root node. As auxiliary r.v.\ we introduce $Q_{n}$, which denotes the number of source-to-sink paths in the series-parallel network corresponding to a forest (i.e., a set) of bucket-recursive trees, where each tree in the forest corresponds to a subblock in the left or right half of the graph. By decomposing the forest into its leftmost tree and the remaining set of trees and taking into account that the number of source-to-sink paths in the forest is the product of the number of source-to-sink paths in the leftmost tree and the corresponding paths in the remaining forest, we obtain the following system of stochastic recurrences:
\begin{equation}\label{eqn:number-path_binary_stochastic}
  P_{n} \stackrel{(d)}{=} Q_{U_{n}}' + Q_{n-2-U_{n}}'', \quad \text{for $n \ge 2$}, \qquad Q_{n} \stackrel{(d)}{=} P_{V_{n}}' \cdot Q_{n-V_{n}}''', \quad \text{for $n \ge 1$}, 
\end{equation}
with $P_{0}=0$, $P_{1}=1$, $Q_{0}=1$, and where the r.v.\ $U_{n}$ and $V_{n}$ are independent of each other and independent of $(P_{n})$, $(P_{n})'$, $(Q_{n})$, $(Q_{n})'$, $(Q_{n})''$ and $(Q_{n})'''$. Furthermore, they are distributed as follows:
\begin{equation*}
  \mathbb{P}\{U_{n}=k\} = \frac{1}{n-1}, \quad 0 \le k \le n-2, \qquad \mathbb{P}\{V_{n}=k\} = \frac{1}{n}, \quad 1 \le k \le n.
\end{equation*}
Introducing $E_{n} := \mathbb{E}(P_{n})$ and $\tilde{E}_{n} := \mathbb{E}(Q_{n})$, the stochastic recurrence above yields the following system of equations for $E_{n}$ and $\tilde{E}_{n}$ (with $E_{0}=0$, $E_{1}=1$ and $\tilde{E}_{0}=1$):
\begin{equation*}
  E_{n} = \frac{2}{n-1} \sum_{k=0}^{n-2} \tilde{E}_{k}, \quad n \ge 2, \qquad \tilde{E}_{n} = \frac{1}{n} \sum_{k=1}^{n} E_{k} \tilde{E}_{n-k}, \quad n \ge 1.
\end{equation*}
Introducing $E(z) := \sum_{n \ge 1} E_{n} z^{n-1}$ and $\tilde{E}(z) := \sum_{n \ge 0} \tilde{E}_{n} z^{n}$ one obtains that $E(z)$ satisfies the following non-linear second order differential equation:
\begin{equation}\label{eqn:DEQ_binary_path-number}
  E''(z) = \frac{1}{1-z} E'(z) + E(z) E'(z), \quad E(0)=1, \quad E'(0)=2.
\end{equation}
Differential equation \eqref{eqn:DEQ_binary_path-number} is not explicitly solvable; furthermore, the so-called Frobenius method to determine a singular expansion fails for $E(z)$. However, it is possible to apply the so-called Psi-series method in the setting introduced in \cite{CheFerHwaMar2014}, i.e., assuming a logarithmic Psi-series expansion of $E(z)$ when $z$ lies near the (unique) dominant singularity $\rho$ on the positive real axis. This yields the following result.
\begin{theorem}
The expectation $\mathbb{E}(P_{n})$ of the number $P_{n}$ of paths from source to sink in a random series-parallel network of size $n$ generated by the binary model has, for $n \to \infty$, the following asymptotic behaviour, with $\rho \approx 0.89\dots$:
\begin{equation*}
  \mathbb{E}(P_{n}) = \frac{2}{\rho^{n}} \cdot \left(1-\frac{\rho^{2}}{(\rho-1)^{2} (n-1) (n-2)} + \mathcal{O}\Big(\frac{\log n}{n^{4}}\Big)\right).
\end{equation*}
\end{theorem}

\end{document}